\documentclass[11pt]{myj}

\usepackage{amsmath}
\usepackage{amsfonts, amssymb, amsthm}

\begin{document}

\title{Counting cycles}
\copyrightinfo{1999}{Igor Rivin}
\author{Igor Rivin}

\address{Mathematics Department, University of Manchester and
  Mathematics Department, Temple University}
\curraddr{UMR 7640 du C.~N.~R.~S., Centre de Math{\'e}matiques, {\'E}cole Polytechnique,
  91128 Palaiseau CEDEX, France} 

\email{rivin@math.polytechnique.fr}

\thanks{I would like to thank {\'E}cole Polytechnique for its support, and
  Ilan Vardi for mentioning the antisemitic question to me, and
  continued encouragement}

\date{\today}


\begin{abstract}
We obtain sharp bounds for the number of $n$-cycles in a finite graph
$G$ as a function of the number of edges, and prove that the complete
graph is optimal in more ways than could be imagined. En route, we
prove sharp estimates on both $\sum_{i=1}^n x_{i}^k$ and $\sum_{i=1}^n
|x_i|^k$, subject to the constraints that $\sum_{i=1}^n x_{i}^2 = C$
and $\sum_{i=1}^n x_i = 0.$
\end{abstract}

\maketitle

\newtheorem{prop}{Proposition}
\newtheorem{defn}{Definition}
\newtheorem{lemma}{Lemma}
\newtheorem{thm}{Theorem}
\newtheorem{cor}{Corollary}
\newtheorem{remark}{Remark}
\newtheorem{ex}{Example}
\newtheorem{question}{Question}
\newtheorem{obs}{Observation}
\newcommand{\tr}{\mathrm{tr\:}}
\newcommand{\spec}{\mathrm{spec\:}}

This note has been inspired by the following question, which had been
asked at the oral entrance exams (see \cite{vardi}) to the Moscow
State University Mathematics Department (MekhMat) to certain applicants:

\begin{question}\label{thequestion}
Let $G$ be a graph with $E$ edges. Let $T$ be the number of triangles
of $G$. Show that there exists a constant $C$, such that for
\textit{all} $G$, $T\leq C E^{3/2}$.
\end{question}

Before proceeding any further, let us answer the question. We will
assume that $G$ is a \textit{simple, loopless, undirected} graph -- that is,
there is exactly one edge connecting two vertices $v$ and $w$ of $G$,
and there are no edges whose two endpoints are actually the same
vertex. 

We will need the following 

\begin{defn}
The adjacency matrix $A(G)$ is a matrix such that
$$(A(G)_{ij}) = \begin{cases}
  1 & \text{if $i$th and $j$th vertices of $G$ are
    adjacent}, \\
    0 & \text{otherwise.}
    \end{cases}
$$
\end{defn}

We shall also need the following observations:

\begin{obs} \label{obs1}
The number of paths of length $k$ between vertices $v_i$ and $v_j$ of
$G$ is equal to $A^k_{ij}$.
\end{obs}

The proof of this is immediate. It follows that the number of
\textit{closed} paths of length $k$ in $G$ is equal to the trace of
$A^k$. Of course, this statement has to be made with some care, since
the trace counts each closed path essentially $2 k$ times (the $2$ is
for the choice of orientation, the $k$ is for the possible starting
points, the ``essentially'' is because this is not true of paths which
consist of the same path repeated $l$ times; that's only counted $2
k/l$ times, or a path followed by retracing the steps backward -- such
a path is counted $k$ times, unless each half is a power of a shorter path,
in which case we count it $k/l$ times...)

Unravelling the various cases, we have:

\begin{equation} \label{eq1}
\tr A = 0,
\end{equation}
since $G$ has no self-loops.
\begin{equation}\label{eq2}
\tr A^2 = 2 E(G),
\end{equation}
where $E(G)$ is the number of edges of $G$.
\begin{equation}\label{eq3}
\tr A^3 = 6 T(G),
\end{equation}
where $T(G)$ is the number of triangles of $G$, and
\begin{equation}\label{eqp}
\tr A^p = 2 p C_p(G),
\end{equation}
where $C_p(G)$ is the number of cycles of length $p$ of $G$ and $p$ is
a prime. For general $k$
\begin{equation}\label{eqk}
\frac{\tr A^k}{k} \leq \text{number of closed paths of length $k$ in
  $G$} \leq \frac{\tr A_k}{2}.
\end{equation}
A much more precise general statement can be made, but this will lead us too
far afield for the moment.

Since $A$ is symmetric, the following observation is self-evident:

\begin{obs}\label{obstr}
$$\tr A^k = \sum_{\lambda \in \spec A} \lambda^k,$$ where $\spec A$ is
the spectrum of $A$ -- the set of all eigenvalues of $A$.
\end{obs}

To answer Question \ref{thequestion} we will also need the following:

\begin{thm}\label{ineqthm}
Let $\mathbf{x}=(x_1, \dots, x_k, \dots)$ be a vector in some Hilbert
space $H$. Then
\begin{equation}\label{holder}
\|\mathbf{x)}\|_p \leq
  \|\mathbf{x}\|_2,
\end{equation}
where $$\|\mathbf{x}\|_q = \left(\sum |x_i|^q\right)^{1/q},$$ and
equality case in the inequality (\ref{holder}) occurs if and only if one
of the $x_i$ is equal to $1$ in absolute value (so that the others are
all equal to $0$.)
\end{thm}

\begin{proof}
It suffices to prove Theorem \ref{ineqthm} under the assumtpion that
$\|\mathbf x\|_2 = 1$ -- the general case follows by rescaling. This
case, however is trivial, and follows from the observation that if
$0\leq y \leq 1$, then $\alpha < \beta$ implies that $y^\alpha \geq
y^\beta$, with equality if and only if $y \in \{0, 1\}$.
\end{proof}

\begin{cor}
Let $M$ be a symmetric matrix. Then
$$(\tr A^k)^2 \leq (\tr A^2)^k,$$ with equality if and only if all the
eigenvalues but one of $A$ are $0$.
\end{cor}

\begin{proof}[Proof of Corollary]
Since, by Observation \ref{obstr}, $\tr A^l = \sum_{\lambda \in \spec
  A} \lambda^l$, this follows immediately from Theorem \ref{ineqthm}.
 \end{proof}

Applying the Corollary in the case $k=3$, together with
eqs. (\ref{eq2}),(\ref{eq3}), we get:

\begin{equation}
(2 E)^{3/2} \geq 6 T,
\end{equation}
 and so
\begin{equation} \label{trestim}
T \leq \frac{2^{1/2}}{3} E^{3/2}.
\end{equation}
We have answered Question \ref{thequestion}, but we have done more: we
found an explicit value for the constant $C$ ($\sqrt{2}/3$), and the
method works without change to show that
\begin{equation}\label{pestim}
C_p \leq \frac{2^{p/2-1}}{p} E^{p/2},
\end{equation}
for prime $p$, while
\begin{equation}\label{kestim}
C_k \leq 2^{k/2-1} E^{k/2}
\end{equation}
in general.

Something not quite satisfying remains about the above argument (aside
from the weak bound for general $k$): it is clear that the equality
case in the estimates (\ref{trestim}) and (\ref{pestim}) is never
attained. This is so, because we know that the equality would
correspond to the spectrum of $A$ consisting of all $0$s and one
non-zero eigenvalue, but this contradicts eq. (\ref{eq1}). So,
potentially we could get a tighter bound by taking (\ref{eq1}) into
account. No easier said than done. We now have the following
optimization problem (for the number of triangles):

Maximize
$$\sum_{i=1}^n \lambda_i^3$$ subject to $$\sum_{i=1}^n \lambda_i = 0,$$
and $$\sum_{i=1}^n \lambda_i^2 = 2 E.$$

This is a typical constrained optimization problem, best approached
with Lagrange multipliers. To avoid (or increase) future confusion, we
note that by scale invariance the optimization problem is equivalent
to

Maximize
$$\sum_{i=1}^n x_i^3$$ subject to $$\sum_{i=1}^n x_i = 0,$$
and $$\sum_{i=1}^n x_i^2 = 1.$$

We know that to find the maximum we need to solve the following
gradient constraint:

$$\nabla(\sum_{i=1}^n x_i^3) = \lambda_1 \nabla(\sum_{i=1}^n x_i) +
\lambda_2 \nabla(\sum_{i=1}^n x_i^2).$$
In coordinates, we have a system of $n$ equations, with the $i$-th
being:
$$E_i:\qquad x_i^2 = \lambda_1 + \lambda_2 x_i.$$
This already tells us that whatever $\lambda_1$ and $\lambda_2$ may
be, there are only two possible values of $x_i$ (independently of
$i$) -- the two roots of the quadratic equation.

Summing all the equations, we get
$$1 = n \lambda_1,$$ so $\lambda_1 = 1/n$.
On the other hand, multiplying $E_i$ by $x_i$ and summing, we see
that:
\begin{equation}\label{lagsum}\sum_{i=1}^n x_i^3 = \lambda_2.
\end{equation}
The left hand side of eq. (\ref{lagsum}) is just the function we are
trying to maximize! It remains, thus, to find a good $\lambda_2.$

Rewriting the equation $E_i$ as
$$x_i^2 - \lambda_2 x_i - 1/n = 0,$$
we obtain:
\begin{equation}\label{quadform}
x_i = \frac{1}{2}\left[\lambda_2 \pm \sqrt{\lambda_2^2 +
    \frac{4}{n}}\right].
\end{equation}
Let us assume that the number of $i$ for which we take the plus sign
in the quadratic formula (\ref{quadform}) exceeds the number of $i$ for
which we take the minus sign by $k$. Summing all of the $x_i$ we get
$$
0 = \sum_{i=1}^n x_i = n \lambda_2 + k \sqrt{\lambda_2^2 +
  \frac{4}{n}}.
$$
(This implies already that $k < 0$.) This translates to the following
equation for $\lambda_2$:

$$\lambda_2^2 = \frac{4 k^2}{n(n^2 - k^2)}.$$ Since we want to make
$\lambda_2$ as large as possible (by eq. (\ref{lagsum})), we want to
make $k^2$ as large as possible on the right hand side. Since at least
one of the $x_i$ has to be negative and at least one positive, $-k$
cannot exceed $n-2$. Thus, the biggest possible value for $\lambda_2$
is
$$\lambda_2 = \frac{n-2}{\sqrt{n (n-1)}},$$
so, after all this work,  we have improved our estimate (\ref{trestim})
to
\begin{equation} \label{trestim2}
T \leq \frac{V-2}{\sqrt{V (V-1)}}\frac{2^{1/2}}{3} E^{3/2},
\end{equation}
($V$ being the number of vertices of our graph $G$).
This is somewhat disappointing: as $E$ (and thus $V$) goes to
infinity, the improvement disappears, and we have the same constant as
before. All the work has not been for nothing, however, for consider
the complete graph on $n$ vertices $K_n$. $E(K_n) = \frac{n (n-1)}{2},$
while $T(K_n) = \frac{n (n-1) (n-2)}{6},$ (since any pair of vertices
defines an edge, while any triple defines a triangle). A simple
computation shows that
\begin{equation}\label{compest}
T(K_n) =  \frac{n-2}{\sqrt{n (n-1)}}\frac{2^{1/2}}{3} E(K_n)^{3/2},
\end{equation}
so the inequality (\ref{trestim2}) is actually an equality in this
case. So the estimate (\ref{trestim2}) is sharp (since it
becomes an equality for an infinite family of graphs), and therefore
constant $\frac{2^{1/2}}{3}$ is also sharp.

A similar calculation (exercise) shows that the constant in the
estimate (\ref{pestim}) is likewise sharp, as shown by the complete
graphs.

A few remarks are in order (as usual).

Firstly, we have inadvertently computed the spectrum of the complete
graph.

The estimate (\ref{trestim2}) and the identity (\ref{compest})
together show that the complete graph $K_n$ is actually maximal
(in terms of the number of triangles) of all graphs with the same
number of vertices and edges as it. This sounds wonderful, until
we realize that it is the \textit{only} graph with $n$ vertices
and $n (n-1)/2$ edges. The identity (\ref{compest}) together with
(\ref{trestim2}) do seem to \textit{suggest} that the complete graph
is maximal (for the number of triangles) of all the graphs with
the same number of edges. We state this as

\begin{question}
Show that the complete graph $K_n$ is the graph containing the most
triangles of the graphs with $\frac{(n-1) n}{2}$ edges.
\end{question}

This question turns out to be not too difficult. The answer is the
content of the following

\begin{thm}\label{compopt}
In a graph $G$ with no more than $n(n-1)/2$ edges, each edge
is contained, on the average, in no more than $n-2$ triangles.
Equality holds only for the complete graph $K_n$.
\end{thm}

\begin{proof}
We will prove the theorem by  induction. Let $v$ be a vertex in
$G$ of maximal degree $d$. Such a vertex is contained in, at most,
$T_v=\min(d(d-1)/2, E(G)-d)$ triangles. This is because there is
at most one triangle per edge connecting two vertices adjacent to
$G$.  and removing it together with the edges incident to it,
leaves a graph $G^\prime$ with $T(G)-T_v$ triangles, $E(G)-d$
edges, and $V(G)-1$ vertices.

Note, first of all, that if the two endpoints of an edge in $G$
have valences $d_1$ and $d_2$, then, if $m = \min(d_1, d_2)$, $E$
is contained in at most $m-1$ triangles. So if the degree of $v$
(assumed to be maximal) was smaller than $n-1$, \textit{no} edge
of $G$ was contained in as many as $n-2$ triangles, so we are
done.

If $d > n-1$, then $G^\prime$ has $n(n-1)/2 -d$ edges, and so each
edge incident to $v$ is contained, on the average, in at most
$\left[n(n-1)-2d\right]/d$ triangles. Now, $$n(n-1)-2d -d(n-2) =
n(n-1) - dn = n(n-1-d) < 0,$$ so the edges incident to $d$ are
contained, on the average, in fewer than $n-2$ triangles. The
number of edges of $G^\prime$ is smaller than $(n-1)(n-2)/2$ (by a
simple calculation), so each of them is contained, on the average,
in at most $n-3$ triangles. Since, at best, each of them was
contained in one more triangle containing $v$, this tells us that
the average was smaller than $n-2$.

If $d=n-1$, repeating the argument as above shows us that for the
 equality to hold $G^\prime$ has to be a complete graph on
$n-1$ vertices, and so $G$ is a complete graph on $n$ vertices.
\end{proof}

Since most numbers are not triangular (triangular numbers being those
of the form $n (n-1)/2$), one can naturally ask the following

\begin{question} Is there a simple characterization of graphs with $k$
  edges which are ``triangle maximal'' (for all $k$)?
\end{question}

and

\begin{question}
Consider all graphs with $E$ edges and $V$ vertices. Is there a way to
characterize the one with the most triangles.
\end{question}

Moving away from graphs as such, the reader will have noted,
perhaps, that our way to maximize the $\sum_{i=1}^n x_i^p$ subject
to the constraints $\|\mathbf x\| = 1$ and $\sum_{i=1}^n x_i = 0$
doesn't work so well for $p \neq 3$, which brings up the
questions:

\begin{question}
Which point $\mathbf{x}$ on the unit sphere $\mathbb{S}^{n-1} \in
\mathbb{R}^n$ and satisfying $\sum_{i=1}^n x_i = 0$ has the
biggest $\sum_{i=1}^n x_i^p$? Which has the biggest $L^p$ norm
(this question is the same of even integer $p$, but quite
different for odd $p$. For non-integer $p$, the first question
doesn't make that much sense...
\end{question}

It turns out that one \textit{can} minimize the sum of $p$-th
powers for $p$ odd. The maximum in this case is attained a the point
satisfying the constraints of largest $L^\infty$ norm.

\begin{thm}\label{oddthm}
The maximal value of $\sum_{i=1}^n x_i^{2p+1}$ subject to the
constraints $\sum_{i=1}^n x_i^2 = 1,$ and $\sum_{i=1}^n x_i = 0$
is attained at the point where $$x_1 = \sqrt{\frac{n-1}{n}}$$ and
$$x_j = -\sqrt{\frac{1}{(n-1)n}} \qquad  j=2,\dots, n.$$ The value
of this maximum is
$$\frac{(n-1)^{2p}-1}{n^{p+1/2}(n-1)^{p-1/2}}.$$
\end{thm}

\begin{proof}
As before, we set up the Lagrange multiplier problem, which has
$n$ equations of the form:
\begin{equation}\label{geneq}
E_i:\qquad x_i^{2p} = \lambda_1 + \lambda_2 x_i.
\end{equation}
Adding all of the equations together, we find that
\begin{equation}\label{l1eq}
n\lambda_1 = \sum_{i=1}^n x_i^{2p},
\end{equation}
while multiplying $E_i$ by $x_i$ and adding the results together
we get
\begin{equation}\label{l2eq}
\lambda_2 = \sum_{i=1}^n x_i^{2p+1},
\end{equation}
so that that sought-after sum is equal to $\lambda_2$, as before.

Further, note that the derivative of $x^{2p}-\lambda_2 x -
\lambda_1$ is equal to $(2p-1) x^{2p-1} - \lambda_2$, which has
exactly $1$ real zero (whatever the value of $\lambda_2$.
Therefore, the equation $x^{2p}-\lambda_2 x - \lambda_1 = 0$ has
at most two real roots. The specifics of our problem are such that
we know that there are exactly two roots, one positive, the other
negative. Call the positive root $\alpha_1$, and the negative root
$\alpha_2$, and suppose that $n_1$ of the $x_i$ are equal to
$\alpha_1$, while $n_2 = n - n_1$ of the $x_i$ are equal to
$\alpha_2$. It follows that
 \begin{equation}\label{aeq}
 \alpha_1 = -
\frac{n_2}{n_1} \alpha_2.
\end{equation}

By eq. (\ref{l1eq}) and (\ref{l2eq}) it follows that $$\lambda_1 =
\frac{1}{n}\left(n_1 \alpha_1^{2p} + n_2
\alpha_2^{2p}\right),\qquad \lambda_2 = n_1 \alpha_1^{2p+1} + n_2
\alpha_2^{2p+1}.$$

 From eq. (\ref{geneq}), we have the following equation for
$\alpha_2$ (where we have substituted for $\alpha_2$ from the
equation (\ref{aeq}):

\begin{equation}
\alpha_2^{2p} = \frac{1}{n}\left(n_1\left(-\frac{n_1}{n_2}
\alpha_2\right)^{2p} + n_2 \alpha_2^{2p}\right)+\left(n_1
\left(-\frac{n_1}{n_2} \alpha_2\right)^{2p+1} + n_2
\alpha_2^{2p+1}\right)\alpha_2.
\end{equation}
Dividing through by $\alpha_2^{2p}$ get

$$1 = \frac{1}{n}\left[\frac{n_2^{2p}}{n_1^{2p-1}} + n_2\right] +
\alpha_2^2\left[-\frac{n_2^{2p+1}}{n_1^{2p}} + n_2\right],$$
 from where, rearranging terms, and replacing $n$ by
$n_1+n_2$, we get

\begin{equation*}
 \alpha_2^2 =
\frac{\frac{1}{n_2} - \frac{1}{n_1+n_2}
\left[\left(\frac{n_2}{n_1}\right)^{2p-1}+1\right]}
{1-\left(\frac{n_2}{n_1}\right)^{2p}} = \frac{n_1}{n_2 (n_1+n_2)},
\end{equation*}
 since, amazingly, everything cancels after clearing
denominators.

So, finally, we see that $$\alpha_2^2 = \frac{n_1}{n_2(n_1+n_2)}$$
while $$\alpha_1^2 = \frac{n_2}{n_1 (n_1+n+2)},$$ thus showing the
first part of the theorem.

 Now, the sum $S$ which we seek is given by
$$S=n_1 \alpha_1 + n_2 \alpha_2 =
\frac{1}{n^{p+1/2}}\left[\frac{n_2^{p+1/2}}{n_1^{p-1/2}} -
\frac{n_1^{p+1/2}}{n_2^{p-1/2}}\right].$$ This is obviously
maximal when $n_2$ is as large as possible, to wit $n-1$, from
which the second part of the theorem follows immediately.
\end{proof}

Notice that since the values of $x_i$ are independent of $p$, it
follows from Theorem \ref{oddthm} that we have proved the following
\begin{thm}
Let $p$ be an odd prime. A graph $G$ with $V$ vertices and $E$ edges has
at most 
$$
\frac{(V-1)^{p-1}-1}{V^{(p+1)/2}(V-1)^{(p-1)/2}}\frac{2^{p/2-1}}{p} E^{p/2}
$$
$p$-cycles, where equality holds if and only if $G$ is the complete
graph $K_{|V|}.$
\end{thm}

\subsection*{Notes on the bibliography} It is hoped that this paper is
reasonably self-contained, however, I would be remiss not to give some
references to related literature. The literature on graph eigenvalues
is vast. For some entry points, the reader is advised to look at the
books of Biggs (\cite{biggs}) and Cvetkovic-Doob-Sachs (\cite{cvet}) for a
general introduction to graph theory, Bollobas' book \cite{bollo} is
good, among many other.

\bibliographystyle{amsplain}

\end{document}